\documentclass[oneside,12pt,reqno]{amsart}

\newtheorem{theorem}{Theorem}[section]

\newtheorem{lemma}[theorem]{Lemma}
\newtheorem{corollary}[theorem]{Corollary}

\theoremstyle{definition}

\newtheorem{definition}[theorem]{Definition}
\newtheorem{example}[theorem]{Example}

\theoremstyle{remark}

\newcommand{\RR}{\mathbb R}
\newcommand{\CC}{\mathbb C}
\newcommand{\QQ}{\mathbb Q}
\newcommand{\ZZ}{\mathbb Z}
\newcommand{\zd}{\mathbb Z^d}
\newcommand{\h}{\mathsf h}
\newcommand{\m}{\mathbf m}
\newcommand{\n}{\mathbf n}
\newcommand{\ve}{\mathbf v}
\newcommand{\be}{\mathbf e}
\newcommand{\bb}{\mathbf b}
\newcommand{\bo}{\mathbf 0}

\begin{document}
\allowdisplaybreaks
\frenchspacing

\bibliographystyle{plain}

\title[Entropy geometry and disjointness]{Entropy geometry and disjointness for
zero-dimensional algebraic actions}

\author{Manfred Einsiedler and Thomas Ward}

\address{ME, University of Washington,
Department of Mathematics, Seattle, WA 98195-4350, USA}
\email{manfred.einsiedler@univie.ac.at}


\address{TW, School of Mathematics, University of East
  Anglia, Norwich NR4 7TJ, United Kingdom}
\email{t.ward@uea.ac.uk}

\thanks{The first
author was supported by the Erwin Schr\"odinger Stipendium J2090 and
the Center for Dynamical Systems at Penn State University}

\begin{abstract}
We show that many algebraic actions of
higher-rank abelian groups on zero-dimensional
groups are mutually disjoint. The proofs
exploit differences in the entropy geometry
arising from subdynamics and a
form of Abramov--Rokhlin formula for half-space entropies.
\end{abstract}

\maketitle

We discuss some mutual disjointness properties
of algebraic actions of higher-rank
abelian groups on zero-dimensional groups.
The tools used are a version
of the half-space entropies introduced by
Kitchens and Schmidt~\cite{MR95f:28022} and adapted
by Einsiedler~\cite{e02},
a basic
geometric entropy formula from~\cite{e02},
and the structure of expansive subdynamics
for algebraic $\mathbb Z^d$-actions
due to Einsiedler, Lind, Miles and Ward~\cite{MR1869066}.
We show that any collection
of algebraic $\mathbb Z^d$-actions on
zero-dimensional groups with entropy
rank or co-rank one that look sufficiently
different are mutually disjoint.
The main results are the following (here $N(\cdot)$
denotes the set of non-expansive directions defined in
Section~\ref{introduction}).

\medskip
\noindent{\bf Theorem~\ref{permtheorem}.}
{\it Let $\mathsf X_1,\dots,\mathsf X_n$ be a collection of
irreducible algebraic zero-dimensional $\mathbb Z^d$-actions, all
with entropy rank one. If
$$
N(\alpha_{j})\backslash\textstyle\bigcup_{k>j}N(\alpha_{k})
\neq\emptyset\mbox{ for }j=1,\dots,n
$$
then the systems are mutually disjoint.}
\medskip

The simplest illustration of Theorem~\ref{permtheorem} is
the fact that
Ledrappier's Example~\ref{ledexample} and its mirror image are
disjoint. This is shown directly in Section~\ref{ledcase} to
illustrate how the Abramov--Rokhlin formula
for half-space entropies may be used.

\medskip
\noindent{\bf Theorem~\ref{main}.}
{\it Let $\mathsf Y$ and $\mathsf Z$ be prime $\zd$-actions with
entropy co-rank one.
If $N(\alpha_Y)\neq N(\alpha_Z)$, then $\mathsf Y$ and
$\mathsf Z$ are disjoint.}
\medskip

Once again the simplest illustration of the meaning
of this result comes from an example of Ledrappier
type: Example~\ref{ledinrank3} is a three-dimensional analogue of
Ledrappier's example.
This is a $\ZZ^3$-action defined by a `four-dot'
condition which has positive entropy $\ZZ^2$-subactions;
it and its mirror image are disjoint.

Surprisingly, it is not
the familiar presence of different non-mixing sets
but the
entropy and subdynamical geometry of the systems that
forces this high level of measurable
difference of structure.
The methods should extend to
entropy rank or co-rank
greater than one,
but the notational and technical difficulties
become more substantial.
Related work for $\mathbb Z^d$-actions by toral
automorphisms has been done by Kalinin and
Katok~\cite{MR1898802}, where more refined information
is found about joinings and the
consequences of the presence of
non-trivial joinings.
Actions by toral automorphisms automatically
have entropy rank not exceeding one.

Our purpose here is to begin to
address some of the problems inherent
in understanding the joinings between
algebraic $\mathbb Z^d$-actions.
The ultimate goal is to extend
results like those of~\cite{MR1898802} to
general algebraic actions,
just as the rigidity results have been
extended from the toral case in~\cite{MR97d:58116}, to
irreducible actions in~\cite{MR2001j:37004}.
In the rigidity theory, entropy rank one
also has a privileged position (see \cite{sb},
\cite{sbks} for the details of how
entropy rank influences rigidity).

Irreducible actions on zero-dimensional groups
are a natural analogue of irreducible actions on finite-dimensional
tori and solenoids, see~\cite{el02}. In particular,
both types of action allow
a local description using locally compact fields.
While $\RR,\CC$ and finite extensions of $\QQ_p$ are used for the
toral and solenoidal cases, for irreducible actions on
zero-dimensional groups locally compact fields
of positive characteristic
are used, namely fields of Laurent series in one variable over a
finite field (see~\cite{MR99b:11089}
and~\cite{el02} for how this works). Using the
local isometry to a product of local fields,
one can define Lyapunov exponents and foliations of the spaces
just as for the toral case. For our purpose it is simpler to
use half-space entropies instead of ultrametric
Lyapunov exponents. Half-space entropies were
introduced in~\cite{MR95f:28022} and adapted to
be defined via state partitions in~\cite{e02}.

A special case showing how the entropy geometry
gives insight into joinings
is dealt with in Section~\ref{ledcase},
and this can be read independently of the rest of the
paper (up to accepting some plausible results on
entropy geometry proved elsewhere).

\section{Introduction}\label{introduction}

An algebraic $\mathbb Z^d$-action is an
action of $\mathbb Z^d$ generated by
$d$ commuting automorphisms of a compact
abelian metrizable group $X$.
Duality (in the sense of Pontryagin)
gives a one--to--one correspondence between
countable modules $M,N,\dots$ over the ring
$R_d=\mathbb Z[u_1^{\pm1},\dots,u_d^{\pm1}]$
and algebraic $\mathbb Z^d$-actions
$\mathsf X_M=(X_M,\alpha_M),\mathsf X_N,\dots$
(see~\cite{MR97c:28041} for an overview
of how this correspondence has been used
to study
algebraic dynamical systems).
It is convenient to write monomials (units)
in $R_d$ in the form
${\mathbf u}^{\mathbf n}=u_1^{n_1}\cdots u_d^{n_d}$.

An algebraic dynamical system $\mathsf X=(X,\alpha)$ automatically preserves
the Haar measure $\lambda=\lambda_X$ on $X$; we
reserve $\lambda$ for Haar measures and $\mu$ for
any $\alpha$-invariant probability measure.

The results on expansive subdynamics we need come
from~\cite{MR1869066}: If $\alpha$ is a $\mathbb Z^d$-action by
homeomorphisms of a compact metric space $(X,\rho)$, then
$N(\alpha)$ denotes the set of {non-expansive} vectors
$\mathbf v\in\mathbb R^d\backslash\{0\}$. That is, $\mathbf v\in
N(\alpha)$ if and only if for every $\epsilon>0$ and $t>0$ there
exists a pair of points $x\neq y$ in $X$ with the property that
$$
\rho\left(\alpha^{\mathbf n}x,\alpha^{\mathbf n}y\right)\le\epsilon
\mbox{ for all }\mathbf n\in
\{\mathbf m\in\mathbb Z^d\mid \mathbf v\cdot\mathbf m<0\}
+t\mathbf v.
$$
The whole action is called {expansive} if there is an
$\epsilon>0$ with the property that
$$
\rho\left(\alpha^{\mathbf n}x,\alpha^{\mathbf n}y\right)\le\epsilon
\mbox{ for all }\mathbf n\in\mathbb Z^d\implies{x=y}.
$$

Let $\alpha$ be an expansive algebraic $\mathbb Z^d$-action on a
zero-dimensional group $X$. By~\cite{e02}, Lemma~2.5, such an
action is automatically an algebraic Markov shift in the following
sense: There are integers $q$ and $s$ and a module of relations
$J\subset (R_d/(q))^s$ such that
\begin{equation}\label{moduleofrelations}
X\cong J^{\perp}\subset
\left(\left(\mathbb Z/q\mathbb Z\right)^s\right)^{\mathbb Z^d},
\end{equation}
where $\cong$ denotes an algebraic isomorphism of $\mathbb
Z^d$-actions and $J^{\perp}$ denotes the annihilator of the
submodule $J$ in the dual group
$\left(\left(\mathbb Z/q\mathbb
Z\right)^s\right)^{\mathbb Z^d}$
of the $R_d$-module $(R_d/(q))^s$.
Under the isomorphism in \eqref{moduleofrelations}, the
$\mathbb Z^d$-action on $X$ corresponds to the natural shift
action on $J^{\perp}$. Having chosen such a presentation of the
system, there is an associated (non-canonical) {state
partition} $\xi=\xi(q,s,J)$ comprising the $q^s$ cylinder sets
obtained by specifying the $\mathbf{0}$ coordinate (some of these
sets may be empty).

Given a $\mathbb Z^d$-action
$\alpha$ by measure-preserving
transformations on $(X,\mu)$ and any measurable partition
$\eta$ of $X$, write
$$
\eta^A=\bigvee_{\mathbf n\in A\cap\mathbb Z^d}\alpha^{-\mathbf n}\eta
$$
for the join of $\eta$ over any set $A\subset\mathbb R^d$. The
{conditional entropy of $A$ given $B$ with respect to $\eta$
and $\mu$} is defined to be $H_{\mu}(\eta^A\vert\eta^B)$. For a
fixed $\eta$ (for instance the state partition for a fixed
presentation), we simply write $H_{\mu}(A\vert B)$ for this
conditional entropy.

The following terminology comes
from~\cite{MR97d:58115} and (in this context)~\cite{MR1869066},
and the resulting condition for
vanishing entropy,
which holds for any invariant measure
$\mu$, is the first key observation in our work.
In the system $\mathsf X_M=(X_M,\alpha_M)$,
a set $A\subset\mathbb R^d$
{codes} $B\subset\mathbb R^d$ if for every $\mathbf m\in
B\cap\mathbb Z^d$ there exists a polynomial
$$
f(\mathbf u)=
\sum_{\mathbf n\in A\cap\mathbb Z^d}
f_{\mathbf n}{\mathbf u}^{\mathbf n}
$$
such that $({\mathbf u}^{\mathbf m}-f)M=0_M$.
Viewing $X_M$ in the form~\eqref{moduleofrelations},
this means that knowledge of the coordinates
$(x_{\mathbf m})_{\mathbf m\in A}$
of a point $x\in X_M$ determines uniquely the coordinates
$(x_{\mathbf m})_{\mathbf m\in B}$.
Notice that
\begin{itemize}
\item $A\mbox{ codes }B\implies H_{\mu}(A\vert B)=0$;
\item $A\mbox{ codes }B\implies A+{\mathbf n}
\mbox{ codes }B+{\mathbf n}$ for every $\mathbf n\in\mathbb Z^d$;
\item $A\mbox{ codes }B, A\cup B\mbox{ codes }C
\implies A\mbox{ codes }B\cup C.$
\end{itemize}

A {joining} of two $\mathbb Z^d$-actions
$\mathsf X_1=(X_1,\mu_1,\alpha_1)\mbox{ and }
\mathsf X_2=(X_2,\mu_2,\alpha_2)$
is a measure $\mu$ on $X_1\times X_2$ invariant
under $\alpha_1\times\alpha_2$ and with the
property that $\mu(A\times X_2)=\mu_1(A)$,
$\mu(X_1\times B)=\mu_2(B)$ for
all measurable $A\subset X_1$,
$B\subset X_2$. Write $J(\mathsf X_1,\mathsf X_2)$
for the collection of all joinings of $\mathsf X_1$ and $\mathsf X_2$.
The systems are {disjoint} if the
only joining is the product measure, so
$J(\mathsf X_1,\mathsf X_2)=\{\mu_1\times\mu_2\}$.

The major simplifying assumption we make is to restrict the
entropy rank: $\alpha$ has {entropy rank one} if there exists
a cyclic subgroup of $\mathbb Z^d$ with positive entropy (viewed
as a $\mathbb Z$-action) but all rank two subgroups of $\mathbb
Z^d$ act with zero entropy. Similarly, $\alpha$ has {entropy
rank} $k<d$ if there is a rank $k$ subgroup of $\mathbb Z^d$
acting with positive entropy (when viewed as a $\mathbb
Z^k$-action) but all subgroups of rank $(k+1)$ act with zero
entropy; finally $\alpha$ has entropy rank $d$ if it has positive
entropy as a $\mathbb Z^d$-action. Similarly, $\alpha$ has entropy
co-rank $k$ if it has entropy rank $(d-k)$. Entropy rank in this
context comes from~\cite{MR1869066}, Sect.~7, and the special
properties of rank one systems are studied in~\cite{el02}
and~\cite{es02}.

A $\mathbb Z^d$-action is called {irreducible} if
it has no closed invariant infinite proper subgroups. Irreducible
actions on connected and zero-dimensional groups are extensively
studied because they exhibit rigidity for $d\ge2$ (cf.~\cite{MR1898802},
\cite{kks}, \cite{MR2001j:37004}). The class of
actions with entropy rank one is a natural extension of the class
of irreducible actions (see~\cite{el02}).

\section{Entropy geometry for $d=2$}\label{entsect}

The results from~\cite{e02} summarized and extended in this
section require the entropy
co-rank to be one. On the other
hand,
many technical simplifications are possible when
the entropy rank is one. In order to have both
conditions, $d=2$ in this section.
We will see in Section~\ref{redstep} that this does
not restrict the applications to rigidity for
larger values of $d$.

\begin{definition}
Let $\mu$ be an invariant measure on the zero-dimensional
expansive algebraic system $\mathsf X=(X,\alpha)$ presented as
in~\eqref{moduleofrelations}. Let $\mathbf v\in\mathbb
R^2\backslash\{0\}$ be a vector with associated half-space
$\mathsf H_{\mathbf v}= \{\mathbf n\in\mathbb Z^2\mid \mathbf
v\cdot\mathbf n<0\}$. The {half-space entropy} of $\mathbf v$
is
\begin{equation}\label{halfspaceentropies}
\mathsf h_{\mu}({\mathbf v})=H_{\mu}(\xi^{\mathbf v^{\perp}}\vert
\xi^{{\mathsf H}_{\mathbf v}})
\end{equation}
where $\xi$ is the state partition (for a fixed presentation) and
$${\mathbf v^{\perp}}=\{\mathbf t\in\mathbb Z^2\mid \mathbf
v\cdot\mathbf t=0\}.$$ If $\mathcal C$ is an $\alpha$-invariant
$\sigma$-algebra, then similarly define the conditional half-space
entropy of $\mathbf v$ to be
$$
\mathsf h_{\mu}(\mathbf v\vert\mathcal C)=
H_{\mu}(\xi^{\mathbf v^{\perp}}\vert
\xi^{{\mathsf H}_{\mathbf v}}\vee\mathcal C).
$$
\end{definition}
For a vector ${\mathbf v}\in\mathbb R^2\backslash\{0\}$, let
$\mathbf v^*$ be a primitive vector in $\mathbb Z^2$ chosen so
that
\[
 \mathsf H_{\mathbf v}+\mathbf v^*=\mathsf H_{\mathbf v}\cup
\mathbf v^{\perp}=\{\mathbf n\in\mathbb Z^2\mid \mathbf
v\cdot\mathbf n\leq 0\}
\]
and let $\ell(\mathbf v,r)$ be chosen so that
$$
\mathsf {\mathbf v}^\perp+(-\ell(\mathbf v,r),\ell(\mathbf
v,r))\mathbf v^*\supseteq \mathbf v^{\perp}+B(r).
$$

The half-space entropy from~\cite{e02} defined
by~\eqref{halfspaceentropies} differs from the entropies used in
\cite{MR95f:28022} in that it depends {\it a priori} on the choice
of presentation~\eqref{moduleofrelations} and only turns out after the
event to be invariant under algebraic isomorphism.
The more robust half-space
entropies in~\cite{MR95f:28022} are automatically invariant under
measurable isomorphism (under suitable hypotheses rigidity makes
measurable and algebraic isomorphism coincide). For Haar measure
the two entropies coincide.

\begin{lemma}\label{hwelldefined} The half-space entropy function
$\mathsf h_{\mu}:\mathbb R^2\backslash\{0\}\to\mathbb R_{\ge0}$ is
independent of the choice of presentation of the
system.\end{lemma}

\begin{proof}
Let $(X,\alpha)$ be an expansive zero-dimensional $\mathbb
Z^2$-action and assume that
$$
X\cong J^{\perp}\subset \left(\left(\mathbb Z/q\mathbb
Z\right)^s\right)^{\mathbb Z^2}
$$
and
$$
X\cong I^{\perp}\subset \left(\left(\mathbb Z/r\mathbb
Z\right)^t\right)^{\mathbb Z^2}
$$
are two presentations of the system giving corresponding state
partitions $\xi$ and $\eta$ with corresponding half-space entropy
functions $\mathsf h_{\mu}^{\xi}$ and $\mathsf h_{\mu}^{\eta}$.
This means that there is an $R_2$-module isomorphism between
$R_2^s/J$ and $R_2^t/I$. Dual to this isomorphism of $R_2$-modules
there is a continuous isomorphism of compact groups from
$I^{\perp}$ to $J^{\perp}$: It follows that there exists an $r>0$
with the property that
$$
\xi^{B(r)}\supseteq\eta\mbox{ and }
\eta^{B(r)}\supseteq\xi
$$
where $B(r)$ is a Euclidean ball of radius $r$ in $\mathbb R^2$
centred at the origin.

Standard properties of entropy and the inclusions
$$\xi^{\mathbf
v^{\perp}}\subset \eta^{\mathbf v^\perp+B(r)} \mbox{ and }
\xi^{\mathsf H_{\mathbf v}}\supset \eta^{\mathsf H_{\mathbf
v}-\ell(\mathbf v, r)\mathbf v^*}
$$
imply that
\begin{equation*}
H_{\mu}\left(\xi^{\mathbf v^{\perp}}\big\vert \xi^{\mathsf
H_{\mathbf v}}\right)\le H_{\mu} \left(\eta^{\mathbf
v^\perp+B(r)}\big\vert \eta^{\mathsf H_{\mathbf v}-\ell(\mathbf
v,r)\mathbf v^*}\right).
\end{equation*}
To obtain a sharper statement, notice that the invariance of the
measure implies (or use~\cite{e02}, Prop.~6.3, for $d=2$)
\begin{eqnarray*}
H_{\mu}\left(\xi^{\mathbf v^{\perp}}\big\vert
\xi^{\mathsf H_{\mathbf v}}\right)&=&\frac{1}{N}H_{\mu}\left(
\xi^{\mathbf v^{\perp}+[0,N)\mathbf v^*}
\big\vert\xi^{\mathsf H_{\mathbf v}}\right)
\\
&\le& \frac{1}{N}H_{\mu}\left( \eta^{\mathbf
v^{\perp}+(-\ell(\mathbf v,r),N+\ell(\mathbf v,r))\mathbf v^*}
\big\vert\eta^{\mathsf H_{\mathbf v}-\ell(\mathbf v,r)\mathbf
v^*}\right)
\\
&\le&\frac{N+2\ell(\mathbf v,r)}{N} H_{\mu}\left(\eta^{\mathbf
v^{\perp}}\big\vert \eta^{\mathsf H_{\mathbf v}}\right).
\end{eqnarray*}
It follows that
$$
\mathsf h_{\mu}^{\xi}(\mathbf v)=
H_{\mu}\left(\xi^{\mathbf v^{\perp}}\big\vert
\xi^{\mathsf H_{\mathbf v}}\right)\le
H_{\mu}\left(\eta^{\mathbf v^{\perp}}\big\vert
\eta^{\mathsf H_{\mathbf v}}\right)=
\mathsf h_{\mu}^{\eta}(\mathbf v),
$$
so by symmetry
$\mathsf h_{\mu}^{\xi}(\mathbf v)=
\mathsf h_{\mu}^{\eta}(\mathbf v).$
\end{proof}
A similar argument shows that the half-space
entropy remains well-defined when conditioned on
an invariant $\sigma$-algebra:
If $\mathcal C$ is a $\sigma$-algebra in $J^{\perp}$
(in the notation of the proof of Lemma~\ref{hwelldefined})
with $\mathcal C^{\prime}$ its image under the isomorphism,
then
\begin{equation}
\h_{\mu}^{\xi}(\mathbf v\vert\mathcal C)= \h_{\mu}^{\eta}(\mathbf
v\vert\mathcal C^{\prime}).
\end{equation}
\begin{example}\label{ledexample}
The archetypal example of a zero-dimensional system with entropy
rank one is due to Ledrappier~\cite{MR80b:28030}: Let
\[
 X_{1}=\bigl\{x\in\mathbb F_2^{\mathbb Z^2}\mid
x_\mathbf n+x_{\mathbf n+\mathbf e_1}+ x_{\mathbf n+\mathbf
e_2}=0\mbox{ for all }\mathbf n\in
 \mathbb Z^2\bigr\},
\]
with $\alpha_{1}$ the $\mathbb Z^2$-action defined by the natural
shift action, and $\lambda=\lambda_{X_1}$ the Haar measure. Then
(cf.~\cite{MR1869066}, Ex.~5.6) $\mathbf v\in N(\alpha_{1})$ if
and only if $\mathbf v$ is parallel to an outward normal of the
convex hull of the set $L=\{(0,0),(0,1),(1,0)\}$. Similarly, the
half-space entropy $\mathsf h_{\lambda}(\mathbf v)$ is positive if
and only if $\mathbf v$ is parallel to an outward normal of the
convex hull of the set $L$.
\end{example}
For a polynomial $f\in R_2$ with $f(\mathbf u)=
\sum_{\mathbf n\in\mathbb Z^2} f_{\mathbf n}{\mathbf u}^{\mathbf n}$, the
{Newton polygon} $\mathcal{N}(f)$ of $f$ is the convex hull of the
{support} $\{\mathbf n\mid f_{\mathbf n}\neq0\}$.

In Example~\ref{ledexample} it is not a coincidence that the set
of points whose convex hull determines the non-expansive
directions is exactly the support of the polynomial $1+u_1+u_2$
generating the module of relations. The same holds more generally
when the entropy co-rank is one -- see~\cite{MR1869066} for the
details.

The following properties hold for any expansive $\mathbb
Z^2$-action $\alpha_M$ on a zero-dimensional
group $X_M$ with entropy rank one, presented as
in~\eqref{moduleofrelations}, and for any $\alpha=\alpha_M$-invariant
measure $\mu$ on $X_M$.
It is useful to talk in terms of {directions}:
a vector $\mathbf v\in\mathbb Z^2\backslash\{0\}$
defines a ray$$r(\mathbf v)=\{t\mathbf v\mid t\in[0,\infty)\};$$
vectors $\mathbf v$ and $\mathbf w$ are in the same
direction if their rays coincide, and a vector $\mathbf v$
is in a rational direction if there is a vector $\mathbf w\in
\mathbb Q^d$ with $r(\mathbf v)=r(\mathbf w)$.
\begin{itemize}
\item There is an {annihilating polynomial}
$f\in R_d$ with the property that $fM=0_M$ and each vertex
coefficient of $f$ is coprime to $q$.
\item For every direction $\mathbf v$, $\mathsf h_{\mu}(\mathbf v)<\infty$.
\item If $\mathbf v$ is not an outward normal vector
to an edge of $\mathcal{N}(f)$, then $\mathsf h_{\mu}(\mathbf
v)=0$.
\item Hence, $\mathsf h_{\mu}(\mathbf v)>0$ only for $\mathbf v$
in finitely many directions, all of them rational.
\end{itemize}

The entropy formula in Theorem~\ref{entformula} relates the
half-space or geometric entropies $\h(\cdot)$ defined
by~\eqref{halfspaceentropies} to the dynamical entropies $h(\cdot)$
of individual elements. In the case of higher entropy rank, an
analogous formula relates the entropy of subactions of the
appropriate rank to geometric entropies of the same rank.

\begin{theorem}\label{entformula}
Let $(X,\alpha)$ be a zero-dimensional algebraic $\mathbb
Z^2$-action with entropy rank one, let $\mu$ be any
$\alpha$-invariant measure on $X$, and let $\mathcal C$ be any
$\alpha$-invariant $\sigma$-algebra. Then
\begin{equation}\label{entformulaeqn}
h_{\mu}(\alpha^{\mathbf n}\vert\mathcal C) =\sum_{\mathbf
v\cdot\mathbf n>0} (\mathbf v\cdot\mathbf n)\mathsf
\h_{\mu}(\mathbf v\vert\mathcal C)
\end{equation}
where the sum is taken over all primitive
integer vectors $\mathbf v$ with $\mathbf v\cdot\mathbf n>0$.
\end{theorem}

The unconditioned version of this is is proved in~\cite{e02};
making the obvious modifications to that proof
shows Theorem~\ref{entformula}. Notice that the
left-hand side is the usual dynamical (conditional) entropy of the
measure-preserving transformation $\alpha^{\mathbf n}$ while the
right-hand side involves only the half-space or geometrical
(conditional) entropies.

The half-space entropies also obey a form of
Abramov--Rokhlin entropy addition formula
(cf.~\cite{MR25:4076}, \cite{MR93m:28023}).

\begin{theorem}\label{t:abramov_rokhlin}
Let $\phi:{\mathsf X}\to{\mathsf Y}$ be a surjective group homomorphism between
zero-dimensional entropy rank one algebraic $\ZZ^2$-systems.
Assume that $\phi$ sends the invariant measure $\mu$ on $X$ to the
invariant measure $\nu$ on $Y$. Then
\begin{equation}\label{arformula}
\mathsf h_{\mu}(\mathbf v)=
\mathsf h_{\nu}(\mathbf v)+\mathsf h_{\mu}(\mathbf v
\vert\phi^{-1}(\mathcal B_Y))
\end{equation}
where $\mathcal B_Y$ denotes the Borel $\sigma$-algebra
on $Y$.
\end{theorem}

\begin{proof}
Assume that $\mathsf X$ and $\mathsf Y$ have been presented in the
form~\eqref{moduleofrelations}, with corresponding state partitions
$\xi$ and $\eta$. In~\eqref{arformula},
$\mathsf h_{\mu}(\cdot)$, $\mathsf h_{\nu}(\cdot)$
are defined using $\xi$, $\eta$ respectively.
Since the half-space entropies are independent of
the chosen presentation of the system we can assume without loss
of generality that $\phi^{-1}\eta\subset\xi$. Then
\begin{eqnarray*}
\mathsf h_{\mu}(\mathbf v)&=&\frac{1}{N}H_{\mu}\left(
\xi^{\mathbf v^{\perp}+[0,N)\mathbf v^*}\big\vert\xi^{\mathsf H_{\mathbf v}}
\right)\\
&=&\frac{1}{N}H_{\mu}\left(
(\phi^{-1}\eta)^{\mathbf v^{\perp}+[0,N)\mathbf v^*}\big\vert
\xi^{\mathsf H_{\mathbf v}}
\right)\\
&&\quad\quad+
\frac{1}{N}H_{\mu}\left(\xi^{\mathbf v^{\perp}+[0,N)\mathbf v^*}\big\vert
\xi^{\mathsf H_{\mathbf v}}
\vee(\phi^{-1}\eta)^{\mathbf v^{\perp}+[0,N)\mathbf v^*}\right)\\
&=&\frac{1}{N}\sum_{n=0}^{N-1}H_{\mu}\left(
(\phi^{-1}\eta)^{\mathbf v^{\perp}+n\mathbf v^*}\big\vert
\xi^{\mathsf H_{\mathbf v}}\vee(\phi^{-1}\eta)^{\mathbf v^{\perp}+
[0,n)\mathbf v^*}\right)\\
&& +\frac{1}{N}\sum_{n=0}^{N-1}H_{\mu}\left(\xi^{\mathbf
v^{\perp}+n\mathbf v^*} \big\vert\xi^{\mathsf H_{\mathbf
v}\cup\left(\mathbf v^{\perp}+ [0,n)\mathbf v^*\right)}\vee
(\phi^{-1}\eta)^{\mathbf v^{\perp}+[0,N)\mathbf v^*}\right)\!.
\end{eqnarray*}
Now
\begin{eqnarray*}
\lefteqn{\frac{1}{N}\sum_{n=0}^{N-1}H_{\mu}\left(
(\phi^{-1}\eta)^{\mathbf v^{\perp}+n\mathbf v^*}\big\vert
\xi^{\mathsf H_{\mathbf v}}\vee
(\phi^{-1}\eta)^{\mathbf v^{\perp}+[0,n)\mathbf v^*}\right)}\\
&&\le\frac{1}{N}\sum_{n=0}^{N-1}H_{\mu}\left(
(\phi^{-1}\eta)^{\mathbf v^{\perp}+n\mathbf v^*}\big\vert
(\phi^{-1}\eta)^{\mathsf H_{\mathbf v}}
\vee(\phi^{-1}\eta)^{\mathbf v^{\perp}+[0,n)\mathbf v^*}\right)\\
&&=H_{\mu}\left((\phi^{-1}(\eta)^{\mathbf v^{\perp}}\big\vert
(\phi^{-1}\eta)^{\mathsf H_{\mathbf v}}\right)\\
&&=H_{\nu}\left(\eta^{\mathbf v^{\perp}}\big\vert
\eta^{\mathsf H_{\mathbf v}}\right)\\
&&=
\mathsf h_{\nu}(\mathbf v).
\end{eqnarray*}
On the other hand, for fixed $n$
$$H_{\mu}\left(\xi^{\mathbf v^{\perp}+n\mathbf v^*}
\big\vert\xi^{\mathsf H_{\mathbf v}\cup\left(\mathbf v^{\perp}+
[0,n)\mathbf v^*\right)}\vee
(\phi^{-1}\eta)^{\mathbf v^{\perp}+[0,N)\mathbf v^*}\right)
$$
$$\rightarrow
H_{\mu}\left(
\xi^{\mathbf v^{\perp}}\big\vert\xi^{\mathsf H_{\mathbf v}}\vee
\phi^{-1}(\mathcal B_Y)\right)$$
as $N\to\infty$ by Martingale convergence.
It follows that
\begin{eqnarray*}
\lefteqn{\frac{1}{N}\sum_{n=0}^{N-1}
H_{\mu}\left(\xi^{\mathbf v^{\perp}+n\mathbf v^*}
\big\vert\xi^{\mathsf H_{\mathbf v}\cup\left(\mathbf v^{\perp}+
[0,n)\mathbf v^*\right)}\vee
(\phi^{-1}\eta)^{\mathbf v^{\perp}+[0,N)\mathbf v^*}\right)}\\
&&\rightarrow
H_{\mu}\left(
\xi^{\mathbf v^{\perp}}\big\vert\xi^{\mathsf H_{\mathbf v}}\vee
\phi^{-1}(\mathcal B_Y)\right)=
\mathsf h_{\mu}(\mathbf v\vert\phi^{-1}(\mathcal B_Y)).
\end{eqnarray*}
This shows that
\begin{equation}\label{myoldman}
\mathsf h_{\mu}(\mathbf v)\le
\mathsf h_{\mu}(\mathbf v)+\mathsf h_{\mu}(\mathbf v
\vert\phi^{-1}(\mathcal B_Y)).
\end{equation}
On the other hand, by the classical Abramov--Rokhlin
entropy addition formula,
\begin{equation}\label{isadustman}
h_{\mu}(\alpha^{\mathbf n})=h_{\mu}(\alpha^{\mathbf n})+
h_{\mu}(\alpha^{\mathbf n}\vert\phi^{-1}(\mathcal B_Y)).
\end{equation}
Equation~\eqref{entformulaeqn} for the trivial
$\sigma$-algebra and the $\sigma$-algebra
$\mathcal{C}=\phi^{-1}\mathcal B_Y$
together with~\eqref{myoldman} and~\eqref{isadustman}
show that
\[
\mathsf h_{\mu}(\mathbf v)= \mathsf h_{\mu}(\mathbf v)+\mathsf
h_{\mu}(\mathbf v \vert\phi^{-1}(\mathcal B_Y)).
\]
\end{proof}

\section{A simple example}\label{ledcase}

In this section we show
how to use the entropy geometry of
Section~\ref{entsect} to prove that Ledrappier's Example~\ref{ledexample},
\[
 X_1=\bigl\{x\in\mathbb F_2^{\mathbb Z^2}\mid x_{\mathbf n}+x_{\mathbf n
+\mathbf e_1}+
x_{\mathbf n+\mathbf e_2}=0\mbox{ for all }\mathbf n\in
 \mathbb Z^2\bigr\},
\]
and its close sibling
\[
 X_2=\bigl\{x\in\mathbb F_2^{\mathbb Z^2}\mid x_{\mathbf n}+x_{\mathbf n
+\mathbf e_1}+
x_{\mathbf n-\mathbf e_2}=0\mbox{ for all }\mathbf n\in
 \mathbb Z^2\bigr\},
\]
are disjoint. That is, if $\alpha_i$ denotes the natural shift
action on $X_i$, and $\mathsf X_i=(X_i,\alpha_i)$, then $J(\mathsf
X_1,\mathsf X_2)=\{\lambda_{X_1}\times\lambda_{X_2}\}$. Let
$\mathsf X=\mathsf X_1\times\mathsf X_2$, and write $\alpha$ for
the Cartesian product of the two $\mathbb Z^2$ shift actions. Let
$\mu$ be a {joining} of the
two systems.

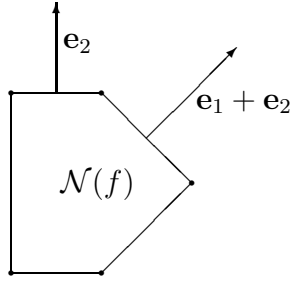
\begin{figure}[ht!]
\setlength{\unitlength}{0.03cm}
\begin{picture}(100,120)
\put(82,73){$\mathbf e_1+\mathbf e_2$}
\put(23,100){$\mathbf e_2$}
\put(22,36){$\mathcal N(f)$}
\put(20,80){\vector(0,1){40}}
\put(60,60){\vector(1,1){40}}
\put(0,0){\circle*{2}}
\put(40,0){\circle*{2}}
\put(80,40){\circle*{2}}
\put(40,80){\circle*{2}}
\put(0,80){\circle*{2}}
\put(0,0){\line(1,0){40}}
\put(40,0){\line(1,1){40}}
\put(0,80){\line(1,0){40}}
\put(40,80){\line(1,-1){40}}
\put(0,0){\line(0,1){80}}
\end{picture}
\caption{\label{fig1}The Newton polygon of the annihilating polynomial}
\end{figure}

A polynomial which annihilates the module corresponding
to $X$ is the product
$$
(1+u_1+u_2)(1+u_1+u_2^{-1})=
u_2^{-1}+u_1u_2^{-1}+u_1^2+u_2+u_1u_2,
$$
with Newton polygon shown in Figure~\ref{fig1}.
Write $\mathcal B_i$ for the Borel $\sigma$-algebra and $\mathcal
N_i$ for the trivial $\sigma$-algebra on $X_i$,
$\xi_i$ for the
state partition in $X_i$ for $i=1,2$, and $\xi=\xi_1\times\xi_2$
for the state partition in $X$.

Part of our purpose here is to show how the half-space
entropies and the
Abramov--Rokhlin formula for half-space entropies
in Theorem~\ref{t:abramov_rokhlin} allow joinings to
be understood.
The first proof below uses the classical
Abramov-Rokhlin formula and the entropy formula
Theorem~\ref{entformula}. The second, much shorter,
proof uses
Theorem~\ref{t:abramov_rokhlin}.

\subsection{Proof of disjointness using
Theorem~\ref{entformula}}
By Section~\ref{entsect},
\begin{equation}\label{eq:sum1}
h_{\mu}(\alpha^{\mathbf e_2})=\mathsf h_{\mu}(\mathbf e_2)
+\mathsf h_{\mu}(\mathbf e_1+\mathbf e_2).
\end{equation}
On the other hand, projecting onto $X_1$ gives a factor of
$\alpha$, so by the Abramov--Rokhlin formula and
Theorem~\ref{entformula}
\begin{eqnarray}
h_{\mu}(\alpha^{\mathbf e_2})&=&h_{\lambda_1}(\alpha_1^{\mathbf
e_2})+h_{\mu}(\alpha^{\mathbf e_2}
\vert\mathcal B_1\times\mathcal N_2)\nonumber\\
&=&h_{\lambda_1}(\alpha_1^{\mathbf e_2})+ \mathsf h_{\mu}(\mathbf
e_2\vert\mathcal B_1\times\mathcal N_2) + \mathsf h_{\mu}(\mathbf
e_1+\mathbf e_2\vert\mathcal B_1\times \mathcal
N_2).\label{eq:sum2}
\end{eqnarray}
Since $\xi_1^{\RR\times(-\infty,0)}=\mathcal
B_1$,
\begin{eqnarray*}
 \mathsf h_{\mu}(\mathbf
 e_2\vert\mathcal B_1\times\mathcal N_2)& =&
 H_\mu(\xi^{\RR\times\{0\}}\vert\xi^{\RR\times(-\infty,0)}\vee
 \mathcal B_1\times\mathcal N_2)\\
 & =&
 H_\mu(\xi^{\RR\times\{0\}}\vert\xi^{\RR\times(-\infty,0)})\\
 &=&\mathsf h_{\mu}(\mathbf e_2).
\end{eqnarray*}
Similarly,
\begin{equation}\label{sumother}
 \mathsf h_{\mu}(\mathbf e_1+\mathbf
e_2\vert\mathcal B_1\times\mathcal N_2)=0,
\end{equation}
and so by comparing~\eqref{eq:sum1}, \eqref{eq:sum2}
and~\eqref{sumother},
\begin{equation}\label{parta}
h_{\lambda_1}(\alpha_1^{\mathbf e_2})= \mathsf h_{\mu}(\mathbf
e_1+\mathbf e_2).
\end{equation}
Projecting onto $X_2$ gives a different factor of $\alpha$
and a similar argument shows that
\begin{equation}\label{partb}
h_{\lambda_2}(\alpha_2^{\mathbf e_2})= \mathsf h_{\mu}(\mathbf
e_2).
\end{equation}
Theorem~\ref{entformula}, \eqref{parta} and \eqref{partb}
together show that
\begin{eqnarray*}
h_{\mu}(\alpha^{\mathbf e_1+\mathbf e_2})&=&\mathsf
h_{\mu}(\mathbf e_2)
+\mathsf h_{\mu}(\mathbf e_1+\mathbf e_2)\\
&=&h_{\lambda_2}(\alpha_2^{\mathbf e_2})+
h_{\lambda_1}(\alpha_1^{\mathbf e_2})\\
&=&\log 4\\
&=&h_{\lambda}(\alpha^{\mathbf e_1+\mathbf e_2}).
\end{eqnarray*}
That is, the joining measure $\mu$ is
a {measure of maximal entropy}
for the transformation $\alpha^{\mathbf e_1+\mathbf e_2}$.
Since $\alpha^{\mathbf e_1+\mathbf e_2}$ is itself an automorphism
of a compact group with finite entropy, it follows
from~\cite{MR40:1582} that $\mu=\lambda=\lambda_{X_1}\times
\lambda_{X_2}$. Thus the
systems $\mathsf X_1$ and $\mathsf X_2$ are disjoint.

\subsection{Proof of disjointness using
Theorem~\ref{t:abramov_rokhlin}}
By the Abramov--Rokhlin formula
for half-space entropies,
\[
 \mathsf h_{\mu}(\mathbf
e_2)=\mathsf h_{\lambda_2}(\mathbf e_2)+\mathsf h_{\mu}(\mathbf
e_2\vert\mathcal N_1\times\mathcal B_2)\geq \log 2,
\]
where we use the fact that $\mathsf h_{\lambda_2}(\mathbf
e_2)=h_{\lambda_2}(\alpha^{\mathbf e_2})=\log 2$. Similarly
\[
 \mathsf h_{\mu}(\mathbf
e_1+\mathbf e_2)=\mathsf h_{\lambda_1}(\mathbf e_1+\mathbf e_2
)+\mathsf h_{\mu}(\mathbf e_1+\mathbf e_2\vert\mathcal
B_1\times\mathcal N_2)\geq \log 2,
\]
so by Theorem~\ref{entformula} the entropy
of the map $\alpha^{\mathbf e_2}$ satisfies
\[
 h_\mu(\alpha^{\mathbf e_2})=\mathsf h_{\mu}(\mathbf e_2)
+\mathsf h_{\mu}(\mathbf e_1+\mathbf e_2)\geq\log
4=h_{\lambda}(\alpha^{\mathbf e_2}).
\]
That is, the joining measure $\mu$ is
maximal for the transformation $\alpha^{\mathbf e_2}$. Since
$\alpha^{\mathbf e_2}$ is itself an automorphism of a compact
group with finite entropy, it follows from~\cite{MR40:1582} that
$\mu=\lambda=\lambda_{X_1}\times \lambda_{X_2}$. Thus the systems
$\mathsf X_1$ and $\mathsf X_2$ are disjoint.

\section{Reduction step}\label{redstep}

In this section we give a corollary to the considerations in
Section~\ref{entsect}, allowing mutual disjointness for entropy
rank one examples to be shown inductively. Recall that an
algebraic $\mathbb Z^d$-action on a zero-dimensional group is
expansive if and only if the corresponding $R_d$-module is
Noetherian (see~\cite{MR97c:28041}). Throughout this section
$\mathsf X$ will be an expansive system.

Recall from~\cite{MR97d:58115} and~\cite{MR1869066}, Sect.~2, the
notion of expansiveness for subsets, and more specifically for
half-spaces $\mathsf H_{\mathbf v}$.
Parameterize half-spaces by the outward normal vector $\mathbf
v$, and write $N(\alpha)$ for the finite set
(see~\cite{MR1869066}, Th.~4.9 and~\cite{el02}, Th.~7.2) of non-expansive
half-spaces.

\begin{theorem}\label{reductionstep}
Let $\mathsf Y=(Y,\alpha_Y,\lambda_Y)$
and $\mathsf Z=(Z,\alpha_Z,\mu_Z)$ be expansive
zero-dimensional algebraic $\mathbb Z^d$-actions with entropy rank
one, and let $\mu$ be in $J(\mathsf Y,\mathsf Z)$.
If there is an integer vector $\mathbf
v\in N(\alpha_Y)\backslash N(\alpha_Z)$, then $\mu$ is invariant
under translation by an infinite subgroup $Y_0\subset Y$. In the
case $d=2$,
$$Y_0=\{y\in Y\mid y_{\mathbf
n}=0\mbox{ for } \mathbf n\in\mathsf H_{\mathbf v}\}.$$
\end{theorem}

Translation in $X=Y\times Z$ by an element $y'\in Y$ means
translation of the form $(y,z)\mapsto(y+y',z)$. Notice that
$\mu_Z$ is any $\alpha_Z$-invariant measure, not
necessarily Haar measure.

\begin{proof}
The first step is to restrict the action to a $\ZZ^2$-subaction
without losing the hypotheses.
By~\cite{MR1869066}, Prop.~7.3, there exists an element
$\alpha^\n$ which acts expansively on $X=Y\times Z$. Let $\m\in\ZZ^d$
be linearly independent to $\n$, and write $P$ for the
plane in $\mathbb R^d$ spanned by $\m$ and $\n$.
Write $\beta$ for the
$\ZZ^2$-subaction generated by $\alpha^{\mathbf k}$
with ${\mathbf k}\in P\cap\zd$.
Similarly, write $\beta_Y$,
$\beta_Z$ for the two factors of $\beta$ on
$Y$ and $Z$. Then $\beta$, $\beta_Y$ and
$\beta_Z$ are each expansive $\ZZ^2$-actions.
We claim the normal vectors to non-expansive
half-spaces for $\beta$ are obtained by
projecting the normal vectors to
non-expansive half-spaces for $\alpha$
onto
the plane $P$ along the orthogonal complement. Thus a half-space in
the plane $P$ is non-expansive if and only if it is contained in a
non-expansive half-space for $\alpha$. This can be seen by a coding argument
similar to the proof of~\cite{MR97d:58115}, Th.~3.6 (replacing
subspaces by half-spaces). Perturbing the plane $P$ slightly does
not affect the expansiveness of the subaction
by~\cite{MR97d:58115}, Lemma~3.4. By
a small perturbation, one can ensure that those
pairs of normal vectors in the finite set $N(\alpha)$ which define
different half-spaces do so in the plane as well. This ensures
that there is a vector $\mathbf v\in N(\beta_Y)\backslash
N(\beta_Z)$. So, without loss of generality assume now that
$\alpha$ is a $\ZZ^2$-action.

Write $\pi_Y:X\to Y$ and $\pi_Z:X\to Z$ for the canonical
projection maps. Then (writing as before $\mathcal B_W$, $\mathcal
N_W$, $\xi_W$ for the Borel $\sigma$-algebra, trivial
$\sigma$-algebra and state partition in $W=Y$ or $W=Z$
respectively)
\begin{eqnarray*}
\mathsf h_{\mu}(\mathbf v)&=&\mathsf h_{\lambda_Y}
(\mathbf v)+\mathsf h_{\mu}(\mathbf v\vert\pi_Y^{-1}(\mathcal B_Y))\\
&=&\mathsf h_{\mu_Z}(\mathbf v)+\mathsf h_{\mu}
(\mathbf v\vert\pi_Z^{-1}(\mathcal B_Z)).
\end{eqnarray*}
Now
$$
\mathsf h_{\mu}(\mathbf v\vert\pi_{Y}^{-1}(\mathcal B_Y))=0
\mbox{ and }
\mathsf h_{\mu_Z}(\mathbf v)=0
$$
since $\mathbf v\notin N(\alpha_Z)$.
It follows that
\begin{equation*}
\mathsf h_{\lambda_Y}(\mathbf v)=\mathsf h_{\mu}(\mathbf v\vert
\pi_Z^{-1}(\mathcal B_Z)),
\end{equation*}
so
\begin{equation}\label{eq:maxequ}
 H_{\lambda_Y}\left(\xi_Y^{\ve^\perp} \vert \xi_Y^{\mathsf H_\ve
 }\right)= H_{\mu}\left(\xi_X^{\ve^\perp} \vert \xi_X^{\mathsf
 H_\ve}\vee\pi_Z^{-1}(\mathcal B_Z)\right){\!\!}.
\end{equation}
We will show that this is the maximal possible value for this
half-space entropy, and deduce the desired translation invariance
property.

Let
$$
Y_0=\{y\in Y\mid y_{\mathbf n}=0\mbox{ for } \mathbf n\in\mathsf
H_\ve\}
$$
and write $\pi:Y_0\to \left(\left(\mathbb Z/q\mathbb
Z\right)^s\right)^{\ve^\perp\cap\ZZ^2}$ for the projection map
onto the coordinates in $\ve^\perp\cap\ZZ^2$ ($Y$ is presented in
the form~\eqref{moduleofrelations} with state partition $\xi_Y$ as
usual). Let
\begin{equation}\label{defetazeta}
\eta_Y=\xi_Y^{\mathsf H_\ve} \mbox{ and } \zeta_Y=\xi_Y^{\mathsf
H_\ve\cup \ve^\perp}{\!}.
\end{equation}
For a measure $\nu$ and partition $\kappa$ write $[x]_{\kappa}$
for the atom of the partition $\kappa$ containing $x$, and
$\nu_{x,\kappa}$ for the associated conditional measure
(characterised by $\int f\mbox{d}\nu_{x,\kappa}= \mathbb
E_\nu(f\vert\kappa)(x)$ for $f\in L^1(\mu)$). By definition of
$\eta_Y$ and $\zeta_Y$ the atom $[y]_{\eta_Y}$ is a union of atoms
$[y+y_0]_{\zeta_Y}$ with $y_0\in Y_0$, where
$[y+y_0]_{\zeta_Y}=[y+y_0']_{\zeta_Y}$ if $\pi(y_0)=\pi(y_0')$.
For the Haar measure $\lambda_Y$ all those $\zeta_Y$-atoms have
the same weight with respect to $\lambda_{y,\eta_Y}$, so that
$$
 H_{\lambda_Y}(\zeta_Y\vert\eta_Y)=
  \log\vert\pi(Y_0)\vert
$$
is finite. The finiteness follows from entropy rank one, see
Section~\ref{entsect}.

We return to the study of $\mu$ on $X=Y\times Z$. Let $\eta_X$ and
$\zeta_X$ be defined similarly to~\eqref{defetazeta},
using the state partition
$\xi_X=\xi_Y\times\xi_Z$. Let
\[
  \eta=\eta_X\vee\pi_Z^{-1}\mathcal B_Z\mbox{ and }
  \zeta=\zeta_X\vee\pi_Z^{-1}\mathcal B_Z.
\]
Then each atom $[x]_\eta$ is a finite union of atoms
$[x+y_0]_\zeta$ with $y_0\in Y_0$, where the sum is
defined by
$x+y=x+(y,0)$. As before, $[x+y_0]_\zeta=[x+y_0]_\zeta$ if
$\pi(y_0)=\pi(y_0')$. By definition, the information function is
$$
I_{\mu}(\zeta\vert\eta)=-\log\mu_{x,\eta}[x]_{\zeta}
$$
and the entropy is its integral
\begin{eqnarray*}
H_{\mu}(\zeta\vert\eta)&=&\int I_{\mu}(\zeta\vert\eta)\mbox{d}\mu\\
&=&\int\sum_{y_0\in\pi(Y_0)}
-\mu_{x,\eta}([x+y_0]_{\zeta})\log\mu_{x,\eta}([x+y_0]_{\zeta})
\mbox{d}\mu.
\end{eqnarray*}
The maximum value of the integral is $\log\vert\pi(Y_0)\vert$,
which is achieved by~\eqref{eq:maxequ}. This happens only
when $\mu_{x,\eta}$ restricted to the partition
$$\{[x+y_0]_\zeta\mid y_0\in Y_0\}$$
of the atom $[x]_\eta$ is a uniform
distribution almost surely. Since translation by $y_0\in Y_0$
permutes the $\zeta$-atoms inside a fixed $\eta$-atom, we deduce
that $\mu(A)=\mu(A+y)$ for any $A\in\zeta$ and $y\in Y_0$. This
argument may be repeated for the next layers, using
$$
\eta'=\eta \mbox{ and } \zeta'=\xi_X^{\mathsf H_\ve+\n},
$$
for some $\n\in\zd\backslash\mathsf H_\ve$. As before a
restricted version of translation invariance for any $A\in\zeta'$
can be shown.
Since this holds for all $\n\in\zd$, it follows that $\mu$ is
invariant under translation by any $y\in Y_0$.
Since $\mathbf v\in\mathsf N(\alpha_Y)$, the subgroup $Y_0$ is
infinite and the theorem follows.
\end{proof}

\section{Applications to disjointness}\label{apps}

The results of Section~\ref{redstep}
suggest the following approach
to mutual disjointness
for systems of this kind. Given a joining
$\mu\in J(\mathsf X_1,\dots,\mathsf X_n)$ of several algebraic systems
$\mathsf X_1,\dots,\mathsf X_n$, look for a vector
$\mathbf v$ that is non-expansive for
$\mathsf X_1$ but expansive for $\mathsf X_2\times\dots\times\mathsf X_n$.
The proof of Theorem~\ref{reductionstep} gives an
equality between two half-space entropies,
and then shows
that $\mu$ is invariant under translation
by a subgroup. If the group is large enough,
this may be enough to deduce that for
almost every $x\in X_2\times\dots\times X_n$,
the conditional measure $\mu_x$ is
Haar measure $\lambda_{X_1}$. This shows
that $\mu=\lambda_{X_1}\times \mu_1$ for
some Borel probability $\mu_1$ on $X_2\times\dots\times
X_n$.
If the process can be repeated with $\mu_1$, then it
shows that the systems are mutually disjoint.

This approach needs two things to happen. First, the non-expansive
sets of the systems must differ enough to keep producing suitable
candidate vectors $\mathbf v$. Second -- a more subtle problem --
the translation invariance provided by Theorem~\ref{reductionstep}
may only give partial information about the measures. To avoid the
latter problem we assume that the systems are irreducible.

\begin{theorem}\label{permtheorem}
Let $\mathsf X_1,\dots,\mathsf X_n$ be a collection of
irreducible algebraic zero-dimensional $\mathbb Z^d$-actions, all
with entropy rank one. If
$$
N(\alpha_{j})\backslash\bigcup_{k>j}N(\alpha_{k})
\neq\emptyset\mbox{ for }j=1,\dots,n
$$
then the systems are mutually disjoint.
\end{theorem}

\begin{proof}
Let $\mu\in J(\mathsf X_1,\dots,\mathsf X_n)$ be a
joining.
Assume by induction that for some $r\ge1$
we know that
$$
\mu=\lambda_{X_{1}}\times\dots\times
\lambda_{X_{r-1}}\times\mu_{r},
$$
and let $\mathbf v$ be a vector in $N(\alpha_{r})\backslash
\bigcup_{k>r}N(\alpha_{k})$. Then $\mu_r\in J(\mathsf
X_r,\dots,\mathsf X_n)$. Apply Theorem~\ref{reductionstep} with
$Y=X_{r}$ and $Z=\prod_{j\neq r}X_{j}$. Since the actions are
irreducible, the subgroup $Y_0$ is dense, so the translation
invariance shows that each fibre of $\mu$ along $Y$ must be Haar
measure on $Y$. That is, $\mu_r=\lambda_{X_r}\times\mu_{r+1}$, and
$$
\mu=\lambda_{X_{1}}\times\dots\times
\lambda_{X_{r}}\times\mu_{r+1},
$$
showing that $\mu=\prod_{j}\lambda_{X_{j}}$ by
induction.
\end{proof}

Since there is a large collection of irreducible polynomials in
$R_2/(p)$ for any fixed prime number $p$,
Theorem~\ref{permtheorem} gives the following
corollary.

\begin{corollary} There is an infinite family
of algebraic $\mathbb Z^2$-actions (on zero-dimensional groups)
with the property that the members of any finite subcollection are
mutually disjoint.
\end{corollary}

\begin{theorem}
Let $\mathsf Y=(Y,\alpha_Y)$ and $\mathsf Z=(Z,\alpha_Z)$
be ergodic expansive
$\mathbb Z^d$-actions with entropy rank one on zero-dimensional
groups. Assume that for any irreducible component of $\mathsf Y$, there is
a vector that is non-expansive on that component, but expansive for
$\alpha_Z$. Then $\mathsf Y$ and $\mathsf Z$ are disjoint.
\end{theorem}

\begin{proof}
 Let $\mathsf X=\mathsf Y\times\mathsf Z$
and let $\mu\in J(\mathsf Y,\mathsf Z)$. Define
 \[
  H_Y=\{y\in Y\mid\mu\mbox{ is invariant under translation by }y\}.
 \]
 If $H_Y=Y$, the measure $\mu$ must be the trivial joining. So
 assume $H_Y\neq Y$ and consider the factors $Y'=Y/H_Y$ and
 $Y'\times Z$ of $Y$ and $X$ respectively;
 the factor measure $\mu'$ is a joining between the Haar measure
 $\lambda_{Y'}$ and $\lambda_Z$. Furthermore, the irreducible
 components of $Y'$ are also irreducible components of $Y$.
 So the assumptions of the theorem remain valid. However,
 by construction the subgroup
 \[
  H_{Y'}=\{y\in Y'\mid\mu'\mbox{ is invariant under translation by }y\}
 \]
 must be trivial. Pick a vector $\ve$ that is non-expansive for
 $\alpha_{Y'}$ but expansive for $\alpha_Z$. By
 Theorem~\ref{reductionstep} the measure $\mu$ is invariant under
translation by an
 infinite subgroup $Y_0\subset Y'$. This contradiction concludes
 the proof.
\end{proof}

\section{Entropy co-rank one in higher dimensions}

In this section we assume that the actions have entropy co-rank
one, allow $d\ge2$, and show disjointness for such actions. The
following replacement for the property of irreducibility is
needed.
Call an algebraic $\mathbb Z^d$-action {prime} if it
is of the form $\mathsf X_M$ for a module $M=R_d/\mathfrak p$ with
$\mathfrak p$ a prime ideal in $R_d$.

\begin{lemma}\label{l:entrk}
 Let $\mathsf Y$ be a prime $\zd$-action with entropy rank
 $k\ge1$. Let $Y'\subset
 Y$ be a closed $\alpha_Y$-invariant subgroup such that the restriction
 $\alpha_{Y'}$ of the action to $Y'$ still has entropy rank $k$.
 Then $Y'=Y$.
\end{lemma}

That is, there are no non-trivial closed invariant subgroups on
which the entropy rank is $k$.

\begin{proof}
This is shown in~\cite{e02}, Proof of Th.~1.2,~Sect.~3.
\end{proof}

\begin{theorem}\label{main}
 Let $\mathsf Y$ and $\mathsf Z$ be prime $\zd$-actions with entropy co-rank
 one.
 If $N(\alpha_Y)\neq N(\alpha_Z)$, then $\mathsf Y$ and
 $\mathsf Z$ are disjoint.
\end{theorem}

In this setting the non-expansive sets are
the set of directions $\mathbf v$ with the property
that the corresponding half-space $\mathsf H_{\mathbf v}$
is non-expansive (see~\cite{MR1869066}, Sect.~2).
In contrast to the case of
entropy rank one, these sets may be infinite.
The next example is the analogue of Section~\ref{ledcase}
for $d=3$.

\begin{example}\label{ledinrank3}
Let
\[
X_1=\bigl\{x\in\mathbb F_2^{\mathbb Z^3}\mid
x_{\mathbf n}+x_{\mathbf n+\mathbf e_1}+
x_{\mathbf n+\mathbf e_2}+
x_{\mathbf n+\mathbf e_3}=0\mbox{ for all }\mathbf n\in
 \mathbb Z^3\bigr\}
\]
and
\[
X_2=\bigl\{x\in\mathbb F_2^{\mathbb Z^3}\mid x_{\mathbf
n}+x_{\mathbf n-\mathbf e_1}+ x_{\mathbf n+\mathbf e_2}+
x_{\mathbf n+\mathbf e_3}=0\mbox{ for all }\mathbf n\in
 \mathbb Z^3\bigr\},
\]
with associated shift $\mathbb Z^3$-actions $\alpha_1$ and
$\alpha_2$. These two systems are associated to
the modules $R_3/(2,f_1)$ and
$R_3/(2,f_2)$ where
$f_1(\mathbf u)=1+u_1+u_2+u_3$ and
$f_2(\mathbf u)=1+u_1^{-1}+u_2+u_3$.
By~\cite{MR1869066} these systems
have entropy co-rank one, and by~\cite{MR97d:58115}, Ex.~2.9,
$N(\alpha_i)$ is the 1-skeleton of the spherical
dual to the Newton polytope $\mathcal{N}(f_i)$ for $i=1,2$.
The vector $\mathbf e_1$ lies in $N(\alpha_2)\backslash N(\alpha_1)$,
so Theorem~\ref{main} shows that $\mathsf X_1$ and $\mathsf X_2$
are disjoint.
\end{example}

The assumption that the entropy co-rank is one
in Theorem~\ref{main} does not seem to be
the whole story, since in Sections~\ref{redstep} and~\ref{apps} we dealt
with general $\zd$-actions with entropy rank one.
Certainly some condition on the entropy rank is required:
If it is allowed to be $d$, then the actions have factors
that are measurably isomorphic to Bernoulli shifts by
\cite{MR96d:22004},
and so have a large space of joinings.
The
geometric picture for entropy rank $k>1$ is more complex. To find a
restriction of the action to a $\mathbb Z^{k+1}$-subactions
without losing the assumptions, one needs a more detailed
description of $N(\alpha)$ -- relating its structure to the
entropy rank of the action -- which is not yet available.

Before we start the proof of Theorem~\ref{main}
we describe the structure of prime actions with entropy co-rank one, and give some
definitions from~\cite{e02}.
If $\mathsf Y$ is
a zero-dimensional prime action with
entropy co-rank one, then
$Y$ is the
dual group of $R_d/(p,f)$ for some prime number $p$ and polynomial
$f$ which is irreducible when considered in $R_d/(p)$ (that the
prime ideal defining the module must have this form when the
entropy co-rank is one follows from~\cite{MR1869066}, Prop.~7.3,
which states that the entropy rank of $\alpha_{R_d/\mathfrak p}$
is equal to the Krull dimension of ${R_d/\mathfrak p}$ if the
characteristic is positive). Clearly $f$ is defined modulo $p$, so
it is natural to assume that $p$ does not divide any nonzero
coefficient of $f$. In the proof of
Theorem~\ref{main} we may assume
that $Z$ is defined in the same way by a prime number
$p'$ and a polynomial $f'$.

Applying a time change for the $\zd$-actions if necessary, we can
make the following simplifying assumptions. Without loss of
generality $-\be_1$ lies in
$N(\alpha_Y)\backslash N(\alpha_Z)$, and
$f\in\ZZ[u_1,u_2^{\pm 1},\ldots,u_d^{\pm 1}]$. The condition
that $-\be_1\in
N(\alpha_Y)$ translates to the property that
$f=f_0+f_1u_1$ for some
$f_1\in\ZZ[u_1,u_2^{\pm 1},\ldots,u_d^{\pm 1}]$ and
$f_0\in\ZZ[u_2^{\pm 1},\ldots,u_d^{\pm 1}]$ which is not a
monomial
by~\cite{MR1869066}, Th.~4.9 and~Ex.~5.7.
Moreover, we can assume that
$f\in\ZZ[u_1,\ldots,u_{d-1},u_d^{\pm 1}]$, and
$f(0,\ldots,0,u_d)=f_0(0,\ldots,0,u_d)\in\ZZ[u_d]$ is not a
monomial
(cf.~\cite{e02}, Lemma~6.8). For $Z$ we can assume that
$f'\in\ZZ[u_1,\ldots,u_d]$, and $f'(0,\ldots,0,u_d)\neq 0$ is a
multiple of a single monomial.

We recall a special case of the notion of lexicographical
half-space entropy for an action of entropy co-rank one
from~\cite{e02}. Let $\Lambda=\ZZ^{d-1}\times\{0\}$ be the
subgroup generated by the first $d-1$ standard basis vectors.
Define lexicographical orders
\begin{eqnarray*}
\m\prec_{\be_d}\n&\mbox{if}&(m_1,\ldots,m_{d-1})
\prec_{\mathrm{lex}}(n_1,\ldots,n_{d-1}),\\
\m\prec\n&\mbox{if}&\m\prec_{\be_d}\n\mbox{ and }m_d=n_d,
\end{eqnarray*}
where $\prec_{\mathrm{lex}}$ is the usual lexicographical order.
Then the lexicographical half-space entropy is defined by
\[
 h_\mu(\be_1,\ldots,\be_{d-1};\be_d)=H_\mu(\xi^{\RR\be_d}|\xi^{S+\RR\be_d}),
\]
where $\xi$ is the state partition and
\[
 S=\{\n\in\ZZ^d\mid\n\succ\mathbf{0}\}\subset\Lambda.
\]
By the above $\h_{\lambda_Y}(\be_1,\ldots,\be_{d-1};\be_d)>0$,
$\h_{\lambda_Z}(\be_1,\ldots,\be_{d-1};\be_d)=0$ --- since
$S+\RR\be_d$ does not code $\RR\be_d$ for $\alpha_Y$, but does for
$\alpha_Z$ and $S+\RR\be_d$ contains $\be_1+H_{-\be_1}$.

Having established these simplifying adjustments and
notations, we turn to the proof of Theorem~\ref{main}.

\begin{proof}
Let $\mu$ be a joining measure, and let $f,f',p,p'$ be as above.
One can change the coefficients of $f$ by multiples of $p$ to
ensure that the non-zero coefficients are all coprime to $pp'$,
and similarly for $f'$. The product $ff'$ annihilates the
$R_d$-module $R_d/(p,f)\oplus R_d/(p',f')$ dual to $X=Y\times Z$,
and every extremal coefficient of $ff'$ is coprime to $pp'$. Thus
$X$ together with $ff'$ and $pp'$ satisfy~\cite{e02}, Lemma~2.5,
and hence the entropy formula in~\cite{e02}, Prop.~6.3; there are
only finitely many directions $\mathbf
w\notin\RR^{d-1}\times\{0\}$ with positive half-space entropies
$h_\mu(\be_1,\ldots,\be_{d-1};\mathbf w)$; the sum
of these half-space entropies
equals the
dynamical entropy $\h_\mu(\alpha_\Lambda)$ of the subaction
defined by $\Lambda$.
Moreover, this entropy formula
remains valid when conditioned by an invariant
$\sigma$-algebra. Just as in Theorem~\ref{t:abramov_rokhlin}, it
follows that
\begin{multline}\label{eq:hrksum1}
 \h_\mu(\be_1,\ldots,\be_{d-1};\be_d)=\\
 \h_{\lambda_Y}(\be_1,\ldots,\be_{d-1};\be_d)+
 \h_\mu(\be_1,\ldots,\be_{d-1};\be_d\vert\mathcal B_Y\times\mathcal N_Z)
\end{multline}
and
\begin{multline}\label{eq:hrksum2}
 \h_\mu(\be_1,\ldots,\be_{d-1};\be_d)=\\
 \h_{\lambda_Z}(\be_1,\ldots,\be_{d-1};\be_d)+
 \h_\mu(\be_1,\ldots,\be_{d-1};\be_d\vert\mathcal N_Y\times\mathcal
 B_Z).
\end{multline}
Coding arguments show that
\begin{eqnarray*}
 \h_{\lambda_Z}(\be_1,\ldots,\be_{d-1};\be_d)&=&0\mbox{ and}\\
 \h_\mu(\be_1,\ldots,\be_{d-1};\be_d\vert\mathcal B_Y\times\mathcal
 N_Z)&=&0.
\end{eqnarray*}
The first equation was noted above; the second follows
similarly. Equations~\eqref{eq:hrksum1} and~\eqref{eq:hrksum2}
imply that
\begin{equation}\label{eq:fineq}
 \h_{\lambda_Y}(\be_1,\ldots,\be_{d-1};\be_d)=
 \h_\mu(\be_1,\ldots,\be_{d-1};\be_d\vert\mathcal N_Y\times\mathcal
 B_Z).
\end{equation}

We use this `maximality property' of the half-space entropy to
deduce a restricted version of translation invariance. Fix
$\ell\geq 1$; let
\begin{eqnarray*}
 U_\ell&=&[0,\ell-1]^{d-1}\times\{0\}\mbox{ and}\\
 S_\ell&=&\ell\{\m\in\Lambda\mid\m\succ\bo\}+U.
\end{eqnarray*}
Define measurable partitions
$\eta=\eta_\ell=\xi^{S_\ell+\RR\be_d}$ and
$\zeta=\zeta_\ell=\xi^{U_\ell+\RR\be_d}\vee\eta$. By
\cite{e02}, Prop.~6.3,
\[
 H_\mu(\zeta\vert\eta\vee\mathcal N_Y\times\mathcal
 B_Z)=\ell^{d-1}\h_\mu(\be_1,\ldots,\be_{d-1};\be_d\vert\mathcal N_Y\times\mathcal
 B_Z)
\]
with a similar expression for the lexicographic half-space entropy
with respect to Haar measure $\lambda_Y$.

Let
\[
 Y_\ell=\{y\in Y\mid y_\n=0\mbox{ for all }\n\in S_\ell+\RR\mathbf
 w\},
\]
and let $\pi:Y_\ell\rightarrow (\ZZ/p\ZZ)^{(U_\ell+\RR\mathbf
w)\cap\zd}$ be the projection map onto the coordinates in
$U_\ell+\RR\mathbf w$. The atom $[x]_{\eta\vee\mathcal
N_Y\times\mathcal B_Z}$ containing the point $x=(x_1,x_2)\in
Y\times Z$ is a subset of $Y\times\{x_2\}$ that splits into many
atoms $[x+y]_{\zeta\vee \mathcal N_Y\times\mathcal B_Z}$ with
$y\in Y_\ell$ (as before $x+y=(x_1+y,x_2)$). Two such atoms for
$y,y'\in Y_\ell$ coincide if and only if $\pi(y)=\pi(y')$, so
there are $|\pi(Y_\ell)|$ such atoms. This gives the upper bound
$\log|\pi(Y_\ell)|$ for the lexicographical half-space entropy,
which is achieved if and only if the conditional measure
$\mu_{x,\eta\vee \mathcal N_Y\times\mathcal B_Z}$ restricted to
the partition
$$\{[x+y]_{\zeta\vee \mathcal N_Y\times\mathcal B_Z}\mid y\in Y_\ell\}$$
is the uniform distribution for $\mu$-a.e. $x\in Y\times Z$. Since
this holds for $\lambda_Y$, the same is true for $\mu$
by~\eqref{eq:fineq}. This implies a restricted translation
invariance property
\begin{equation}\label{transinvstatement}
 \mu(A+y)=\mu(A)\mbox{ for } A\in\eta_\ell\mbox{ and }y\in Y_\ell.
\end{equation}

Let
\[
 Q_\ell=\{\n\in\zd\mid n_i\geq -\ell\mbox{ for all }i<d \},
\]
and $\m=\ell\be_1+\cdots+\ell\be_{d-1}$. Then $Q_\ell+\m\subset
(S\cup U)+\RR\be_d$, and so
$\alpha^{-\m}\xi^{Q_\ell}\subset\eta_\ell$.

Let
\begin{eqnarray*}
 T&=&\{\n\in\zd\mid n_i\geq 0\mbox{ for some }i<d\}\mbox{ and}\\
 Y_0&=&\{y\in Y\mid y_\n=0\mbox{ for all }\n\in T\}.
\end{eqnarray*}
As above $\alpha^{-\m}Y_0\subset Y_\ell$. Therefore
$\alpha$-invariance of the measure allows us to reformulate
\eqref{transinvstatement} as
\begin{equation}\label{evenmoretransinvstatement}
 \mu(A+y)=\mu(A)\mbox{ for }A\in\xi^{Q_\ell}\mbox{ and }y\in Y_0.
\end{equation}
However, $\bigcup_\ell Q_\ell=\zd$, and
so~\eqref{evenmoretransinvstatement} implies that
$\mu(A+y)=\mu(A)$ for $y\in Y_0$ and every measurable $A\subset
Y\times Z$.

To complete the proof of the theorem, we need to show that $\mu$
is in fact invariant under translation by all $y\in Y$. Let
$Y'\subset Y$ be the closure of the group generated by the orbit
of $Y_0$ under the action, and let $\alpha_{Y'}$ be the
restriction of the action to the invariant subgroup $Y'\subset Y$.
The invariance of $\mu$ under $\alpha$ and under translation by
$Y_0$ implies that $\mu$ is invariant under translation by $Y'$.
We claim that the subaction $(\alpha_{Y'})_\Lambda$ has positive
entropy; this shows that $\alpha_{Y'}$ has entropy rank $d-1$, and
Lemma~\ref{l:entrk} shows that $Y'=Y$.

Suppose $Y_0$ is the trivial subgroup. Then the restriction map
\[
 \varphi:Y\rightarrow (\ZZ/(p))^T
\]
to the coordinates in $T$ is injective (that is, the dual groups
$R_d/(p,f)$ and $\ZZ[\mathbf u^\n\mid\n\in T]/(p,f)$ are equal).
Therefore for $\m=-\be_1-\cdots-\be_{d-1}$ there exists a
polynomial $g\in\ZZ[\mathbf u^\n\mid\n\in T]/(p,f)$ with
$$\mathbf u^\m-g\in (p,f).$$
We will show that this contradicts the special geometry of $f$ and
$T$. In the following the equations are meant modulo $p$, so
suppose $\mathbf u^\m-g=hf$ for some polynomial $h$. Split $h$
into a sum $h=h'+h''$ with $h''\in\ZZ[\mathbf u^\n\mid\n\in T]/(p)$
and $h'\in\ZZ[\mathbf u^\n\mid n_i<0$ for all $i<d]$. Taking the
product and using $f\in\ZZ[u_1,\ldots,u_d]$ gives $hf=h'f+h''f$
and $h''f\in\ZZ[\mathbf u^\n\mid\n\in T]/(p)$. Since $g\in\ZZ[\mathbf
u^\n\mid\n\in T]/(p,f)$, we must have $h'f\in\mathbf u^\m+\ZZ[\mathbf
u^\n\mid\n\in T]/(p)$. Let $h'_{\operatorname{min}}$ be the sum of
those coefficients of $h'$ whose exponent $\n$ of $\mathbf u$ is
minimal with respect to $\prec_{\be_d}$. Let
$f_{\operatorname{min}}$ be the same for $f$. Then by the
assumption on $f$ the polynomial
$f_{\operatorname{min}}\in\ZZ[u_d^{\pm 1}]$ cannot be a single
monomial. The terms of $h'f$ whose exponents are minimal are
exactly the terms in
$h'_{\operatorname{min}}f^{\phantom'}_{\operatorname{min}}$. Since
the latter is contained in $\ZZ[\mathbf u^\n\mid n_i'<0$ for all
$i<d]$, it must be equal to $\mathbf u^\m$, which is a
contradiction since $f_{\operatorname{min}}$ is not a monomial.

By the above $Y_0\subset Y'$ is nontrivial, which implies that
$$\h_{\lambda_{Y'}}(\bb_1',\ldots,\bb_{d-1}';\mathbf w)>0$$ and so
by the entropy formula~\cite{e02}, Prop.~6.3,
$h_{\lambda_{Y'}}\bigl((\alpha_{Y'})_\Lambda\bigr)>0$ as claimed.
\end{proof}


\end{document}